\documentclass{elsart}

\usepackage{amsmath}
\usepackage{amssymb}
\usepackage{cite}
\usepackage{array}
\usepackage{graphicx}
\usepackage{enumerate}
\usepackage[center,small]{caption2}
\theoremstyle{plain} 

\theoremstyle{plain} 

\theoremstyle{plain} 

\theoremstyle{definition} 

\theoremstyle{definition} 
\journal{Arxiv}
\begin{document}

\begin{frontmatter}



\title{\textbf{An Evolution Strategy Method for Optimizing Machining Parameters of Milling Operations}}


\author{Zhixin Yang}
\address{Department of Mathematics, University of Wisconsin-Eau Claire,\\ Eau Claire, WI 54702, USA }

\begin{abstract}
In this paper, an evolutionary strategy (ES) method is introduced as
an optimization approach to solve problems in the manufacturing
area. The ES method is applied to a case study for milling
operations. The results show that it can effectively yield good
results.
\end{abstract}

\end{frontmatter}

\newpage
\framebox{Insert Nomenclature here}
\section{Introduction}

In manufacturing industry, it is important to determine the optimal
machining parameters in order to maximize total profit rate and to
increase the quality of the final product for machining operations.
Traditionally, the selection of parameters is carried out by the
experience of planners or with the help of machining data handbooks.
These may not guarantee the optimum performance and minimization of
costs. Therefore, a number of researchers have tried to deal with
the optimization of machining parameters using different approaches.
Compared with some deterministic methods
\cite{petropoulos1973305,armarego1994437,gupta1995115,wang200177}
evolutionary algorithms (EAs) are more attractive to engineers since
EAs are robust, effective and easy to implement,

Many different types of EAs have been proposed for optimizing
machining parameters of milling operations, such as genetic
algorithm (GA) \cite{yildiz20062041}, simulated annealing (SA)
\cite{baskar20051078}, ant colony algorithm \cite{baskar20051078},
immune algorithm \cite{yildiz2009224}, and some hybrid algorithms
\cite{yildiz2012asoc,yuan2008924,yuan2010640,yuan201036}.

Evolution Strategies(ESs), originally developed by Rechenberg and
Schwefel\cite{Rec,Sch}, are algorithms which imitate the principle
of natural evolution as a method to solve parameter optimization
problems\cite{ZM,yuan2008257}. ESs can reject infeasible individuals
directly (also called ``death penalty''). It is probably the easiest
way to handle constraints. Since in the model introduced in sequel,
the feasible search space constitutes a reasonably large portion of
the whole search space. This strategy seems feasible. Moreover, it
is also computationally efficient, because when a certain solution
violates a constraint, it is assigned a fitness of zero. Therefore,
no further calculations are necessary to estimate the degree of
infeasibility of such a solution. To the best knowledge of the
authors, no one has used ESs to optimize machining parameters of
milling operations. So in this paper, we use an ES to solve this
problem. The results show that it can effectively yield good
results.

\section{Mathematical model}
Depth of cut, feed rate and cutting speed have the greatest effect
on the success of the machining operation. Depth of cut is usually
predetermined by the workpiece geometry and operation sequence. It
is recommended to machine the workpiece with the required depth in
one pass to keep machining time and cost low, when possible.
Therefore the problem of determining the machining parameter is
reduced to select the proper cutting speed and feed rate
combination. The mathematical model developed by M. Tolouei-Rad et
al.\cite{tolouei19971} is considered in this work. The model is also
considered in several papers (\cite{baskar20051078, yildiz2012asoc,
yildiz2009261, yildiz2012}).

\subsection{Objective function}
The main focus of this work is to maximize the total profit rate and
can be determined by
\begin{equation}\label{eq1}
  P_r=\frac{S_p-C_u}{T_u}.
\end{equation}
The unit cost can be represented by
\begin{equation}
\begin{split}
  C_u=& c_{mat}+(c_l+c_o)t_s+\sum^m_{i=1}(c_l+c_o)K_{1i}V_i^{-1}f^{-1}_i\\
      & +\sum^m_{i=1}c_{ti}K_{3i}V_i^{(1/n)-1}f_i^{[(w+g)/n]-1}+\sum^m_{i=1}(c_l+c_o).
\end{split}
\end{equation}
The unit time to produce a part in the case of multi-tool milling
can be defined by
\begin{equation}
  T_u=t_s+\sum^m_{i=1}K_{1i}V^{-1}_if^{-1}_i+\sum^m_{i=1}t_{tci}
\end{equation}
\subsection{Constraints}
In practice, possible range of cutting speed and feed rate are
limited by the following constraints
\begin{enumerate}[1.]
\item Maximum machine power
\item Surface finish requirement
\item Maximum
cutting force permitted by the rigidity of the tool
\item Available
feed rate and spindle speed on the machine tool
\end{enumerate}
\subsubsection{Power}
The machining parameters should be selected such that maximum
machine power is used. The required machining power should not
exceed available motor power. Therefore the power constraint can be
written as
\begin{equation}
  C_5Vf^{0.8}\leq 1,
\end{equation}
where
\begin{equation}
  C_5=\frac{0.78K_pWz a_{rad}a}{60\pi deP_m}.
\end{equation}
\subsubsection{Surface finish}
The required surface finish $R_a$, must not exceed the maximum
attainable surface finish $R_{a(at)}$ under the conditions.
Therefore the surface finish for end milling becomes
\begin{equation}
  C_6f\leq 1,
\end{equation}
where
\begin{equation}
  C_6=\frac{318[\tan(la)+\cot(ca)]^{-1}}{R_{a(at)}}.
\end{equation}
And for end milling
\begin{equation}
  C_7f^2\leq 1,
\end{equation}
where
\begin{equation}
  C_7=\frac{318(4d)^{-1}}{R_{a(at)}}.
\end{equation}
\subsubsection{Cutting force}
The total cutting force $F_c$ resulting from the machining operation
must not exceed the permitted cutting force $F_c$(per) that the tool
can withstand. The permitted cutting force for each tool has been
considered as its maximum limit for cutting forces. Therefore the
cutting force constraints becomes
\begin{equation}
C_8F_c\leq 1,
\end{equation}
where
\begin{equation}
C_8=1/F_c\rm{(per)}.
\end{equation}
\subsubsection{Speed limits}
\begin{enumerate}[1.]
\item Face milling: 60--120 m/min
\item Corner milling: 40--70m/min
\item Pocket milling: 40--70 m/min
\item Slot milling1: 30--50 m/min
\item Slot milling2: 30--50 m/min
\end{enumerate}
\subsubsection{Feed rate limits}
\begin{enumerate}[1.]
\item Face milling: 0.05--0.4mm/tooth
\item Corner milling: 0.05--0.5mm/tooth
\item Pocket milling: 0.05--0.5mm/tooth
\item Slot milling1: 0.05--0.5mm/tooth
\item Slot milling2: 0.05--0.5mm/tooth
\end{enumerate}
\section{The ES method}
\subsection{Introduction of ESs}
Evolution Strategies can be understood as `intelligent'
probabilistic search algorithms which are based on the evolutionary
process of biological organisms in nature. 


The procedure of one type of ES can be described as follows. For the
$\mu$ initial individuals, in optimization terms, each individual in
the population is encoded by a real number which represents a
possible solution to a given problem. Then a recombination procedure
and mutation procedure are executed to produce new `offspring'(i.e.
{\em children}) solutions with $\eta(\eta>\mu)$ individuals. The
fitness of each individual in `offspring' solutions is evaluated
with respect to a given objective function. After evaluation
procedure the offspring solutions are sorted and the last $\eta-\mu$
individuals are deleted. This reproduction-evaluation-selection
cycle is repeated until a satisfactory solution is found.

A more comprehensive overview of ESs can be found, e.g. in
\cite{ZM,B&S} and references therein. The applications of ES could
be found in
\cite{yuanarxiv2014,yuan201311408,zhu20131,yuan2010640,yuan2008924,yuan201036,yuan2008257}

The basic steps of the procedure can be shown as:
\begin{quote}
\texttt{Generate an initial population;}\\
\rm{\textbf{repeat}}:\\
\texttt{\mbox{\hspace{.5cm}}Recombinant and mutate individuals to produce children; \\
\mbox{\hspace{.5cm}}Evaluate fitness of the children;\\
\mbox{\hspace{.5cm}}Select the population from the children;} \\
\rm{\textbf{until}} \texttt{a satisfactory solution has been found.}
\end{quote}

By means of the above ES procedure, a computing method
for optimizing machining parameters of milling operations is
described in the following parts of the section.

\subsection{Representation and fitness function}
We denote $X$ as a vector composed of two arrays representing
cutting speed and feed rate, i.e., $x_1=V$ and $x_2=f$,
$X=(x_1,x_2)$. The target function \eqref{eq1} is taken as the
fitness function in the ES.


\subsection{Initial population, recombination and mutation}
\subsubsection{Initial population and Recombination}
The initial population consists of $\mu$ individuals
$\overrightarrow{X}(0)=\{X^{(1)}(0),\cdots,X^{(\mu)}(0)\}$, where
$X^{(i)}(0)=(X^i, \sigma^i)$. The initial components of $X$ are
generated randomly from feasible domains. All initial $\sigma_i$ are
valued by $3.0$ here.

Select two individuals:
\begin{align*}
(X^1,\boldsymbol{\sigma}^1)&=((x_1^1,\cdots,x_l^1),(\sigma_1^1,\cdots,\sigma_l^1))\;\\
\text{and}\ \ \
(X^2,\boldsymbol{\sigma}^2)&=((x_1^2,\cdots,x_l^2),(\sigma_1^2,\cdots,\sigma_l^2)).
\end{align*}
There are two types of recombination operators:
\begin{itemize}
\item {\textbf{discrete}}, where the new offspring is
\[
(X,\boldsymbol{\sigma})=((x_1^{q_1},\cdots,x_l^{q_l}),(\sigma_1^{q'_1},\cdots,\sigma_l^{q'_l}))
\]
with $q_i$ and $q'_i$ equal either 1 or 2;
\item {\textbf{intermediate}}, where the new offspring is
\[ (X,\boldsymbol{\sigma}) =((\alpha
x_1^1+(1-\alpha)x_1^2,\cdots,(\alpha x_l^1+(1-\alpha)x_l^2),
(\alpha\sigma_1^1+(1-\alpha)\sigma_1^2,\cdots,
\alpha\sigma_l^1+(1-\alpha)\sigma_l^2)), \]
with $\alpha\in(0,1)$.
\end{itemize}
By the suggestion of Schwefel\cite{Sch}, the discrete recombination
operator is executed on $X$, and the intermediate recombination
operator is executed on $\boldsymbol{\sigma}$.

\subsubsection{Mutation}
Apply mutation to the offspring $(X,\boldsymbol{\sigma})$, the
resulting new offspring $(X',\boldsymbol{\sigma}')$ is obtained,
where
\begin{align}\label{43}
\sigma_i'&=\sigma_i\cdot\exp(\tau'\cdot
N(0,1)+\tau\cdot N_i(0,1))\;\notag\\
\text{and} \ \ \ x'_i&=x_i+\sigma'_i\cdot N_i(0,1),
\end{align}
in which, \\
$N(0,1)$ ----- a random number reconciled standard normal distribution;\\
$N_i(0,1)$ ----- a random number reconciled standard normal
distribution aimed at $i$th component; \\
$\tau'$ ----- global coefficient;\\
$\tau$ ----- local coefficient;\\
$l$ ----- the number of components in $X$.

According to the suggestion of Schwefel\cite{Sch}, $\tau'$ is valued
by $(\sqrt{2l})^{-1}$ and $\tau$ is valued by
$(\sqrt{2\sqrt{l}})^{-1}$ in the computation.

Notice that $X$ has box constrains. So the offspring generated by
mutation procedures may not be feasible.
Then we just let it take the bound value, i.e. if
$x_i>\overline{x_i}$ (or $< \underline{x_i}$), let
$x_i=\overline{x_i}$ (or $=\underline{x_i}$).

The recombination and mutation steps will not stop until $\eta$
offsprings are generated.

\subsection{Evaluation, Selection and Stop criterion }
We evaluate every $(X,\boldsymbol{\sigma})$ by its fitness function
and sort them. If an individual is infeasible, the fitness value is
0. Then we choose first $\mu$ individuals as new parents. In our ES
method, we let $\mu=15$, $\eta=105$ and set a variable to record the
current best fitness value. If a fitness value is better than the
record, we will update the record. The algorithm will stop if the
record keeps unchanged after 1000 iterative loops.

\section{Case Study}
The component as shown in Figure 1 is to be produced using a CNC
milling machine. It is desirable to find the optimum machining
parameters, which result in the maximum profit rate. Specifications
of the machine, material and values of constants are given below.
Also, the geometric information on the required operations and tools
is presented in Tables 1 and 2.

\framebox{Insert Figure 1 here}

\textbf{Constants:}
\begin{quote}
$S_p = \$25$\\
$cmat =\$0.50$\\
$c_o = \$1.45$ per min\\
$c_l = \$0.45$ per min\\
$t_s = 2$ min\\
$t_{tc} = 0.5$min\\
$C = 33.98$ for HSS tools\\
$C = 100.05$ for carbide tools\\
$w = 0.28$\\
$K_p = 2.24$\\
$W = 1.1$\\
$n = 0.15$ for HSS tools\\
$n = 0.3$ for carbide tool\\
$g = 0.14$
\end{quote}
\textbf{Machine tool data:} Type: Vertical CNC milling machine $P_m
= 8.5$kW, $e = 95\%$\\
\textbf{Material data:} Quality: 10L50 leaded steel. Hardness = 225
BHN

\framebox{Insert Table 1 and Table 2 here.}

We use our ES to solve this case and compare it to some other
methods in the literature. The results listed in Table 3 show that
ES can obtain good results similar to hybrid methods.

\framebox{Insert Table 3 here.}

\section{Conclusion}

In this paper, we used an ES method to optimize machining parameters
of milling operations. The results showed that it can effectively
give good results and it can be a good alternative in similar
problems in engineering.

\newpage

\begin{figure}
\caption{An example}
  \includegraphics[width=80mm]{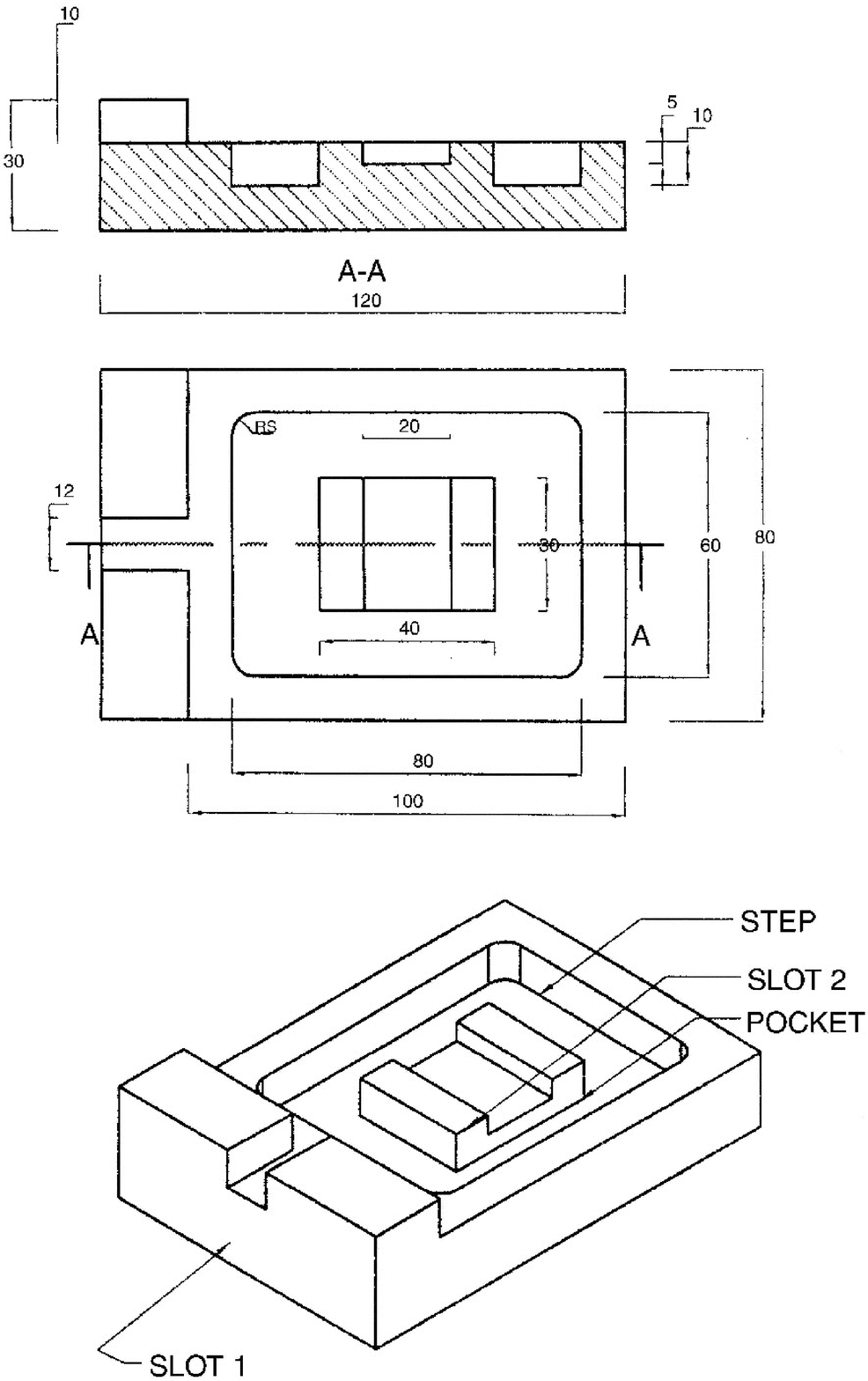}
\end{figure}

\begin{table}
\begin{tabular}{ll}
\hline
\textbf{Nomenclature} &\\
$a$, $a_{rad}$ & Axial depth of cut, radial depth of cut (mm)\\
$ca$ & Clearance angle of the tool (degrees)\\
$C_i$ $(i=1,\ldots,8)$ & Coefficients carrying constants values\\
$c_l$, $c_o$ & Labour cost, overhead cost (\$/min)\\
$c_m$, $c_{mat}$, $c_t$ & Machining cost, cost of raw material per
part, cost
of a cutting tool (\$)\\
$C_u$ & Unit cost (\$)\\
$d$ & Cutter diameter (mm)\\
$e$ & Machine tool efficiency factor\\
$F$ & Feed rate (mm/min)\\
$f$ & Feed rate, (mm/tooth)\\
$F_c$, $F_c$(per) & Cutting force, Permitted cutting force (N)\\
$G$, $g$ & Slenderness ratio, exponent of slenderness ratio.\\
$K$& Distance to be travelled by the tool to perform the
operation (mm)\\
$K_i(i=1,2,3)$ & Coefficients carrying constant values\\
$K_p$ & Power constant depending on the workpiece material\\
$la$ & Lead (corner) angle of the tool (degree)\\
$m$ & Number of machining operations required to produce
the product \\
$N$ & Spindle speed (rev/min)\\
$n$ & Tool life exponent\\
$P$, $P_m$ & Required power for the operation, motor power
(kW)\\
$P_r$ & Total profit rate (\$/min)\\
$R$& Sale price of the product excluding material, setup
and tool changing costs (\$)\\
$R_a$, $R_{a(at)}$ & Arithmetic value of surface finish, and
attainable
surface finish ($\mu$m)\\
$S_p$ & Sale price of the product (\$)\\
$T$, $T_u$ & Tool life (min), Unit time (min)\\
$t_m$, $t_s$, $t_{tc}$ & Machining time, set-up time, tool changing
time (min)\\
$V$, $V_{hb}$, $V_{opt}$ & Cutting speed, recommended by handbook,
optimum (m/min)\\
$w$ & Exponent of chip cross-sectional area\\
$W$ & Tool wear factor\\
$z$ & Number of cutting teeth of the tool\\
\hline\\
\end{tabular}
\end{table}

\begin{table}
\caption{Required machining operation}
  \begin{tabular}{cccccc}
    \hline
    Operation Number & Operation type & Tool number & $a$(mm)& $K$(mm) &
    $R_a$($\mu$m)\\
    \hline
    1 & Face milling & 1 & 10 & 450& 2\\
    2 & Corner milling &2& 5 &90 & 6\\
    3 & Pocket milling &2& 10 &450& 5\\
    4 & Slot milling &3& 10& 32& --\\
    5 & Slot milling &3& 5& 84& 1\\
    \hline
  \end{tabular}

\end{table}

\begin{table}
\caption{Tools data}
  \begin{tabular}{cccccccc}
    \hline
    Tool Number & Tool type & Quality & $d$(mm)& $z$ &
    Price(\$)& $la$ & $ca$\\
    \hline
    1 & Face mill & Carbide & 50 & 6& 49.50 &45& 5\\
    2 & End mill &HSS& 10 &4 & 7.55& 0& 5\\
    3 & End mill &HSS& 12 &4& 7.55 & 0& 5\\
    \hline
  \end{tabular}

\end{table}
\begin{table}
  \caption{Comparison of the results for milling operation}
  \begin{tabular}{cccc}
   \hline
   Method & $C_u$-Unit cost & $T_u$-Unit time & $P_r$-Profit Rate\\
   \hline
   Handbook \cite{handbook} &\$18.36 & 9.40 min & 0.71/min\\
   Method of feasible direction \cite{tolouei19971} &\$11.35 & 5.48 min &
   2.49/min\\
   Genetic algorithm \cite{yildiz2009261} & \$11.11 & 5.22 min & 2.65/min\\
   Ant colony algorithm \cite{baskar20051078} & \$10.20 & 5.43 min &
   2.72/min\\
   Hybrid particle swarm \cite{yildiz2012} & \$10.90 & 5.05 min &
   2.79/min\\
   Immune algorithm & \$11.08 & 5.07 min & 2.75/min\\
   Hybrid immune algorithm \cite{yildiz2009261} & \$10.91 & 5.07 min & 2.79/min\\
   Hybrid differential evolution algorithm \cite{yildiz2012asoc} & \$10.90 & 5.00 min & 2.82/min\\
   Evolutionary strategy & \$10.91 & 5.00 min & 2.82/min\\
   \hline
  \end{tabular}
\end{table}
\end{document}